\documentclass[12pt]{article}

\usepackage{amsmath}
\usepackage{amssymb}
\usepackage{amsfonts}

\usepackage{amsthm}

\newtheorem{theorem}{Theorem}
\newtheorem{lemma}{Lemma}

\newcommand{\lbl}{\label}

\newcommand{\cM}{{\cal M}}

\def\N{\mathbb{N}}
\def\Z{\mathbb{Z}}

\def\R{\mathbb{R}}
\def\E{\mathbb{E}}

\def\Pr{\mathbb{P}}

\def\0{{\bf 0}}

\def\bS{{\bf S}}
\def\bZ{{\bf Z}}

\renewcommand{\E}{\mathbb E \,}

\newcommand{\tod}{\stackrel{{\cal D}}{\longrightarrow}}

\newcommand{\toP}{\stackrel{{P}}{\longrightarrow}}

\newcommand{\tovia}{\rightsquigarrow}
\renewcommand{\leftrightarrow}{\leftrightsquigarrow}

\newcommand{\Var}{{\rm Var}}

\newcommand{\tG}{{\tilde{G}}}
\newcommand{\tPr}{{\tilde{\P}}}
\newcommand{\tE}{{\tilde{\E}}}

\newcommand{\tp}{{\tilde{p}}}

\newcommand{\surv}{\texttt{surv}}

\newcommand{\eps}{\varepsilon}
\def\bdm{\begin{displaymath}}
\newcommand{\edm}{\end{displaymath}}
\def\benu{\begin{enumerate}}
\def\eenu{\end{enumerate}}
\def\beqn{\begin{equation}}
\def\eeqn{\end{equation}}
\def\be{\begin{equation}}
\def\ee{\end{equation}}
\def\bea{\begin{eqnarray}}
\def\eea{\end{eqnarray}}
\newcommand{\bean}{\begin{eqnarray*}}
\newcommand{\eean}{\end{eqnarray*}}
\newcommand{\bear}{\begin{eqnarray}}
\newcommand{\eear}{\end{eqnarray}}

\renewcommand{\epsilon}{\varepsilon}

\def\R{\mathbb{R}}

\renewcommand{\P}{{\mathbb P}}

\def\qed{\hfill\hbox{${\vcenter{\vbox{
    \hrule height 0.4pt\hbox{\vrule width 0.4pt height 6pt
    \kern5pt\vrule width 0.4pt}\hrule height 0.4pt}}}$}}

\begin{document}
\title{\bf The strong giant in a
 random digraph}

\author{Mathew D. Penrose$^{1}$\\
{\normalsize{\em University of Bath}} }

\maketitle


 \footnotetext{ $~^1$ Postal address: Department of
Mathematical Sciences, University of Bath, Bath BA2 7AY, United
Kingdom}
\footnotetext{ Email address: {\texttt m.d.penrose@bath.ac.uk} }




\begin{abstract}
Consider a random directed graph on $n$ vertices with
 independent identically distributed
 outdegrees with distribution $F$ having mean $\mu$,
 and destinations of arcs selected uniformly at random.
We show that if $\mu >1$ then for large $n$ there is
very likely to be a unique giant strong component with proportionate
size given as the product of two branching
process survival probabilities, one with offspring distribution
$F$ and the other with Poisson offspring distribution
with mean $\mu$.  If $\mu \leq 1$ there is very likely to be no
giant strong component. We also extend this to
allow for $F$ varying with $n$. 
\end{abstract}


\section{Introduction}

Given
$n \in \N := \{1,2,3,\ldots \}$, and
 given a probability distribution $F$ on $\Z_+ := \N \cup \{0\}$,
 consider a random directed multigraph
 $G_{n,F}$ on vertex set
 $[n] := \{1,\ldots,n\}$, defined as follows
(a multigraph is a graph with multiple arcs and loops allowed).
Let each vertex of $G_{n,F}$ have outdegree independently
sampled from distribution $F$.
Given these outdegrees, the arcs
have their destinations sampled independently
uniformly from $[n]$.

Consider also a random directed graph $\tilde{G}_{n,F}$, similar
to $G_{n,F}$ but with loops and multiple arcs excluded.
Let the outdegrees be sampled independently
from distribution $F$ as before, 
 denoting the
respective outdegrees by $\xi_1,\ldots,\xi_n$. 
Given these outdegrees, let the set of destinations of the
arcs from vertex $i$ in $\tG_{n,F}$ be selected
uniformly at random from the collection of
all $\binom{n-1}{\xi_i}$ subsets of
$[n]\setminus \{i\}$ with $\xi_i$ elements, independently
of the arcs from other vertices. If $\xi_i > n-1$,
 include all arcs from $i$ in the graph $\tG_{n,F}$
(so in this case the outdegree of $i$ is $n-1$, not $\xi_i$).

In the special case where $F$ is
 the Dirac distribution at $k$ for
some $k \in \N$ (i.e. $F(\{k\}) =1$), the 
random directed graph $\tG_{n,F}$ is also known as
$G_{n,k-{\rm out}}$ or just $G_{k-{\rm out}}$, mentioned
in \cite{Boll} and studied in 
\cite{FF} and elsewhere.

Random graph models with specified degree distributions
(for example, with power-law decay of the tails)
are of much recent interest, and directed graphs are
often a better model for real-world networks than the
undirected ones (see \cite{Broder,NSW,DMS} and references therein). 
Our model is a simple and natural way to allow for an
arbitrary specified outdegree distribution in a random directed graph.

We use the abbreviation `digraph' to mean either
a directed graph or a directed multigraph.
For vertices $i,j$ of a finite digraph $G$ we write
$i \tovia j$ if  
there is a directed path from $i$ to $j$ (or if $i =j$), and
$i \leftrightsquigarrow j$ if both 
$i \tovia j$ and $j \tovia i$. 
We say $G$ is {\em strongly connected}
if for any two vertices $i,j$  we have
$i \leftrightsquigarrow j$.
For $k \in \N$
let $L_k(G)$ denote the number of vertices in
the $k$th largest 
strongly connected component
of $G$ (if $k$ exceeds the number of such components,
set $L_k(G)=0$.) 

Set $\mu_F := \sum_k k F(\{k\})$, 
the mean of distribution $F$.
In Theorem \ref{thgiant}, we show that $G_{n,F}$
and $\tG_{n,F}$ enjoy
a `giant component' phenomenon also seen 
 in other random graph models such as the
Erd\"{o}s-R\'enyi
 random graph (see e.g. \cite{Boll}), if and only if $\mu_F>1$.
In Theorems \ref{thFvary} and \ref{thunif} 
we extend this result to allow for 
$F$ varying with $n$, and in 
Theorem \ref{CDSthm} we give a further result
on the distributional limit of the proportionate
number of vertices $j \in [n]$ such that $1 \tovia j$.

 Related random digraph models are considered
(and results analogous to Theorem \ref{thgiant} are derived)
 in \cite{BGJ}  and \cite{CF}, but they
are not the same as ours.
In \cite {CF} the degrees are imposed
globally whereas here they are determined locally.
In \cite{BGJ}, each vertex has  a randomly determined
type, and each arc is included at random with probability
determined by the type of its endpoints.

 The random digraph model considered here
seems at least as natural as those in
\cite{BGJ} and \cite{CF}. 
In some sense, it is intermediate between
the one in \cite{CF} which is homogeneous
and the one in \cite{BGJ} which is inhomogeneous; loosely speaking,
one may say that a random graph is homogeneous if all of
its vertices have the same status (see \cite{BGJ} and references
therein).
Our graphs are {\em semi-homogeneous}
in the sense that they are inhomogeneous with respect to 
outdegree but homogeneous with respect to indegree. 
Nevertheless, our model does not appear to be a special
case of that in \cite{BGJ}, since the status (present/absent)
of
different arcs from a given vertex are not conditionally independent
given  the type (i.e. outdegree) of that vertex.  Also,
in \cite{BGJ}  the number of types is
assumed to be finite, whereas we allow for $F$ with infinite support.

\section{Statement of results}
\label{secstate}

Given a probability distribution $F$ on $\Z_+ \cup \{\infty\} $,
let $x_F$ be the smallest solution
in $[0,1]$ of the equation $x = \phi_F( x)$,
where we set $\phi_F(x) = \sum_{k=0}^\infty x^k F(\{k\})$, 
and set $\sigma(F) := 1-x_F$.
It is well-known (see e.g. \cite{Feller}) that
$\sigma(F)$ is the survival probability of
a Galton-Watson branching process with offspring distribution
$F$, and that
$\sigma (F) > 0$
if and only if $\mu_F > 1$.
  In the special case where $F$ is a Poisson
distribution with parameter $\mu \in [0,\infty)$
(so $\phi_F(x) = e^{\mu (x-1)}$), we write
$\sigma'(\mu)$ for $\sigma(F)$, and we set $\sigma'(\infty) =1$.

Let $\toP$ denote convergence in probability.

\begin{theorem}
\label{thgiant}
Given any probability distribution $F$ on $\Z_+$,
as $n \to \infty$ we have 
\bea
n^{-1} L_1(G_{n,F})  \toP  \sigma'(\mu_F) \sigma(F);   ~~~~~
 n^{-1}  L_2(G_{n,F})  \toP 0,
\label{M6}
\eea
and 
\bea
n^{-1} L_1(\tG_{n,F})  \toP  \sigma'(\mu_F) \sigma(F);
~~~
 n^{-1} L_2(\tG_{n,F})  \toP 0. 
\label{L2conv}
\eea
\end{theorem}
It is natural to ask whether the convergence 
in probability statements of Theorem \ref{thgiant}
hold uniformly over all choices of the outdegree
distribution $F$.
This amounts to asking whether similar statements hold
if we allow $F$ to vary with $n$.
Our next results tell us that this is indeed the case
if for each $n$ we impose a deterministic bound $b_n$ on the 
outdegrees in $G_{n,F}$, satisfying $b_n = o(n)$ (i.e. $b_n/n \to 0$
as $n \to \infty$).
For $n \in \N$, let $\cM_n$
be the class of probability distributions $F$
on $\Z_+$ which are supported by $\{0,1,\ldots,n\}$, i.e.
which satisfy $F(\{0,1,\ldots,n\}) =1$.

Given a probability distribution $F$ on $\Z_+ \cup \{\infty\}$, and
a sequence of probability distributions $(F_n)_{n \geq 1}$ on
$\Z_+$, we write
$F_n \Longrightarrow F$ if
$F_n$ converges weakly to $F$ as $n \to \infty$, i.e.
if $\lim_{n \to \infty} F_n(\{k\}) = F(\{k\})$ for all $k \in \Z_+$. 
We note that if $F_n \Longrightarrow F$ then
$\sigma(F_n) \to \sigma (F)$ as $n \to \infty$.
Likewise, $\sigma'(\mu)$ is continuous in $\mu$, including
at $\mu = \infty$.

\begin{theorem}
\label{thFvary}
Let $(b_n)_{n \geq 1}$ be an $\N$-valued sequence
with $b_n = o(n)$ as $n \to \infty$.
Suppose $(F_n)_{n \geq 1}$ is a sequence of probability distributions
on $\Z_+$ with $F_n \in \cM_{b_n}$
for each $n$, satisfying $F_n \Longrightarrow F$ for some
 probability distribution $F$ on $\Z_+ \cup \{\infty\}$, and
also $ \mu_{F_n} \to \mu_\infty$ as $n \to \infty$ for some
 $\mu_\infty \in [0,\infty]$.
Then
as $n \to \infty$ we have 
\bea
n^{-1} L_1(G_{n,F_n})  \toP  \sigma'(\mu_\infty) \sigma(F);
~~~~~~
n^{-1}  L_2(G_{n,F_n})  \toP 0
,
\label{M15}
\eea
and 
\bea
n^{-1} L_1(\tG_{n,F_n})  \toP  \sigma'(\mu_\infty) \sigma(F)
;
~~~~~~
n^{-1}  L_2(\tG_{n,F_n})  \toP 0
.
\label{M15a}
\eea
\end{theorem}
As a corollary, we may deduce a result about uniform
convergence. To state  this we need a metrization
of convergence in probability. Given random variables
$X,Y$ on the same probability space, set
$d(X,Y) : = \sup \{ \eps : \Pr [ |X-Y| > \eps ] > \eps \}$. 

\begin{theorem}
\label{thunif}
Let $(b_n)_{n \geq 1}$ be an $\N$-valued sequence
with $b_n = o(n)$ as $n \to \infty$.
Then
\bea
\lim_{n \to \infty} \sup_{F \in \cM_{b_n}}
  d (n^{-1} L_1(G_{n,F}) ,    \sigma'(\mu_F) \sigma(F) ) =0,
\label{M16}
\\
\lim_{n \to \infty} \sup_{F \in \cM_{b_n}}
  d (n^{-1} L_1(\tG_{n,F}) ,    \sigma'(\mu_F) \sigma(F) ) =0,
\label{M16a}
\eea
and
$\lim_{n \to \infty} \sup_{F \in \cM_{b_n}}
  n^{-1}  (d(L_2(G_{n,F}) , 0)) =
\lim_{n \to \infty} \sup_{F \in \cM_{b_n}}
  n^{-1}  (d(L_2(\tG_{n,F}) , 0)) =0$.
\end{theorem}
Theorem \ref{thgiant} cannot be deduced directly
from Theorem \ref{thFvary} or \ref{thunif}, because
the distribution $F$ in the statement of Theorem \ref{thgiant}
could have unbounded support.
We shall give an example at the end of this section to show
 that 
that the condition $b_n = o(n)$ is needed in
Theorems \ref{thFvary} and \ref{thunif};
it is too much to expect the
convergence in Theorem \ref{thunif} to be uniform
over {\em all} probability distributions on $[n]$.

Our proof of Theorems \ref{thgiant}, \ref{thFvary} and \ref{thunif} relies
heavily on the following result, which is of independent
interest. Given a vertex $i$ of a finite digraph $G$,
let $T_i:= T_i(G)$ be the number of vertices $j$ such that
$i \tovia j$ (including $i$ itself). 
In epidemic modelling, this may be considered 
as the final size of the epidemic, i.e. the
total number of individuals  to become 
infected starting from a single infected individual $i$.  
Let
$\tod$ denote convergence in distribution. 
\begin{theorem}
\lbl{CDSthm}
Suppose $(F_n)_{n \geq 1}$ is a sequence
of probability distributions on $\N_0$ 
such that $F_n \Longrightarrow F$ for
some
 probability
distribution $F$ on $\N_0 \cup \{\infty\}$ and 
$\mu_{F_n}  \to \mu_\infty$ for some $\mu_\infty  \in 
[0,\infty]$.
Suppose either that $F_n =F$ for all $n$, or that
there exists an $\N$-valued sequence $(b_n)_{n \in \N}$ 
such that $b_n = o(n)$ as $n \to \infty$ and
$F_n \in  \cM_{b_n}$ for all $n$. Then
\bea
n^{-1} T_1 ( G_{n,F_n}) \tod \sigma'(\mu_\infty) \xi
\label{CDSeq}
\eea
and
\bea
n^{-1} T_1 ( \tG_{n,F_n}) \tod \sigma'(\mu_\infty) \xi,
\label{CDS2eq}
\eea
where $\xi$ is a Bernoulli random variable with parameter
$\sigma (F)$.
\end{theorem}
Theorem \ref{CDSthm} extends a recent result of
Comets, Delarue and Schott  \cite{CDS}, who
proved (\ref{CDSeq})   
 in the case with $F_n=F$ for all $n$.
Only the case of (\ref{CDSeq}) with $\mu_\infty < \infty$ 
(but with $F_n$ possibly varying with $n$) is
used in proving our other results; the rest of Theorem
\ref{CDSthm} is included for its own sake.

In the rest of this paper we  prove the theorems stated above.
Before embarking on the detailed proof
we introduce further notation
and give some intuitive ideas  behind the proof of the theorems.

Given $F$, let 
$(Z_m)_{m \geq 0} 
:= (Z_m(F))_{m \geq 0}$ 
be a Galton-Watson
 branching process with offspring distribution $F$.
If $F(\{\infty\}) >0$ then we may have $Z_m = \infty$ for
some $m$, in which case we set $Z_n = \infty $ for all $n \geq m$.

Given also $\mu_\infty \in [0,\infty]$,
 let $(Z'_m)_{m \geq 0}
: = (Z'_m(\mu_\infty))_{m \geq 0}
 $  be a branching process with Poisson offspring
 distribution with mean $\mu_\infty$, independent
of $(Z_m)_{m \geq 0}$ (with $Z_0= Z'_0=1$). 
Let $T := \sum_{m=0}^\infty Z_m$ and $T' := \sum_{m=0}^\infty Z'_m$.
If $\mu_\infty = \infty$  then set $Z'_m := +\infty$ for all $m \geq 1$, 
 and set $T' = \infty$.
Then   $ \sigma(F) = \Pr[T = \infty]$
and
 $ \sigma'(\mu_\infty) = \Pr[T' = \infty]$.

Given a vertex $i$ of a digraph $G$,
for $m \in \N$ let 
$S_{i,m}(G)$ denote the size of the $m$th  out-generation starting from $i$,
i.e. the number of vertices $j$  of $G$ such that there
is a directed path from $i$ to $j$ and the shortest
such path is of length $m$. 
Let
$S'_{i,m}(G)$ denote the size of the $m$th  in-generation starting from $i$,
i.e. the number of vertices $j$  of $G$ such that there
is a directed path from $j$ to $i$ and the shortest
such path is of length $m$. 
Set $S_{i,0}(G) =1$ and $S'_{i,0}(G)=1$. 
Set ${\bf S}_{i,m}(G) := (S_{i,0}(G),\ldots,S_{i,m}(G))$
and
set ${\bf S}'_{i,m}(G) := (S'_{i,0}(G),\ldots,S'_{i,m}(G))$.
 Then $T_i(G) = \sum_{m=0}^\infty S_{i,m}(G)$.
 Set $T'_{i} := T'_i(G) := \sum_{m=0}^\infty S'_{i,m}$, the
total number of vertices that can be reached by
a backwards directed path from vertex $i$ in the graph $G$. 
Let $\Pr_{n,F}$ (respectively 
$\tPr_{n,F}$)
denote probability
with reference to the graph $G_{n,F}$ (respectively $\tG_{n,F}$).
Let $\E_{n,F}$ (respectively $\tE_{n,F}$) 
denote expectation
with reference to the graph $G_{n,F}$ (respectively $\tG_{n,F}$).

The intuition for Theorem \ref{thgiant}
or \ref{thFvary}
is that (with $G_{n,F}$ or $\tG_{n,F}$ or $G_{n,F_n}$ or
$\tG_{n,F_n}$  now denoted $G_n$ for short)
for any fixed $m$ and $i$,
the distribution of the random vector  ${\bf S}_{i,m}(G_n) $
 approximates in the large-$n$ limit to
 that of the branching process $\bZ_m := (Z_0,\ldots,Z_m)$.
Moreover the indegree of vertex $i$ in $G_n$
 is asymptotically   
Poisson with mean $\mu_\infty$ (where in the setting
of Theorem \ref{thgiant} we set $\mu_\infty = \mu_F$), and the random vector
 $\bS'_{i,m} (G_n) $ 
 converges in distribution
(as $n \to \infty$)
 to  the random vector 
$\bZ' := (Z'_0,\ldots,Z'_m)$, 
with $\bS_{i,m}(G_n) ,\bS'_{i,m}(G_n) ,\bS_{j,m}(G_n) ,$
and $\bS'_{j,m}(G_n) $
asymptotically independent for fixed $i,j,m$ with
 $j \neq i$.
We justify these assertions in Lemma \ref{lembranch} below.

One might reasonably hope that for large $K$,
the condition that $T_{i} (G_n) > K$
and $T'_{i}(G_n)  > K$
would be approximately necessary and sufficient 
for $i$ to lie in a giant strong component. 
Our argument to
demonstrate this (in Lemma \ref{CDSlem} below)  is based
on the branching process approximation combined with 
Theorem \ref{CDSthm}.


We now give an example to show what can go wrong
if we drop the condition $b_n= o(n)$ in
Theorem \ref{thFvary}, \ref{thunif} or \ref{CDSthm}.
Suppose that we take $F_n(\{n-1\}) = 2n^{-1}$
and $F_n(\{2\} ) = 1- 2n^{-1}$.
Then the limiting distribution $F$ of $F_n$
is a unit point mass at 2 (with $\sigma(F)=1$) 
and the limit of $\mu_F$ is 4.
If the conclusion of Theorem \ref{thFvary}
were still true for this example, then the $n^{-1}L_1(\tilde{G}_{n,F_n})$
should approximate to $\sigma'(4)$.

Consider, however, the successive in-generations $S'_{1,m}(\tilde{G}_{n,F_n}), m \geq
 1$.
While
 $S'_{1,1}(\tilde{G}_{n,F_n})$,
does converge
in distribution to the first  generation $Z'_1$
of a branching process with Poisson offspring distribution with mean
4,
 the second generation
 $S'_{1,2}(\tilde{G}_{n,F_n})$ does not converge in distribution
to $Z'_2$. This is because the  vertices of  
 $S'_{1,1}(\tilde{G}_{n,F_n})$ decompose into two
types, namely those of outdegree 2 and those of outdegree
$n-1$ (with an asymptotically Poisson number of each type with mean 2),
 but in subsequent generations 
 $S'_{1,m}(\tilde{G}_n)$ for $m \geq 2$, there are no
vertices of the second type (because all such vertices
would be included in the first generation). Therefore
the branching process approximation fails.

In our proofs we shall repeatedly use the fact
 that for any $n \in \N$ and  any probability distribution
$F$ on $\Z_+$, the random digraph $\tG_{n,F}$ stochastically
dominates $G_{n,F}$, i.e.
\bea
G_{n,F} \prec_{\rm st} \tG_{n,F},  
\label{eqmonot}
\eea
 in the sense that  there exist
 coupled realizations of these two random digraphs
 for which $G_{n,F}$ (with loops removed and
multiple edges reduced to single edges) 
is a (directed) subgraph of $\tG_{n,F}$.

\section{Proof of Theorem \ref{CDSthm}}
\label{sec.CDSpf}

Throughout this section we assume that
  $F_n$, $F$, $\mu_\infty$ and
(if applicable) $b_n$ are as in the statement of 
Theorem \ref{CDSthm}.
Also we write just $\Pr_n$ (respectively $\tPr_n$, $\E_n$,
$\tE_n$) for
 $\Pr_{n,F_n}$ (respectively $\tPr_{n,F_n}$,   $\E_{n,F_n}$,
 $\tE_{n,F_n}$).

Given a digraph $G= (V,E)$, and given $i,j \in V$ we
write $i \to j$ if there is an arc of $G$ from
$i$ to $j$. Given also $B \subset V$ we write
$i \to B$ if $i \to j$ for at least one $j \in B$.
In the following lemma the notation $H$ stands for
`hit' and $A$ stands for `avoid'.

\begin{lemma}
\label{lemPoconv}
Fix $r,s \in \Z_+ $,
let $H_{r,s}$ be the event
$\{ 1 \to \{ r+2,\ldots, r+1 +s \}  \}$ and
  let $A_r$ 
be the complement of the event $H_{0,r}$.
Then
\bea
\lim_{n \to \infty} \Pr_{n,F_n} [A_r] = \lim_{n \to \infty} \tPr_{n,F_n}
[A_r] =1,
\label{limAvoid}
\eea
and
\bea
\lim_{n \to \infty} (n \Pr_{n,F_n}[ H_{r,s}
| A_r ]) =
\lim_{n \to \infty} (n \tilde{\Pr}_{n,F_n}
[ H_{r,s }  | A_r ]) =
s  \mu_\infty.
\label{limHit}
\eea
\end{lemma}
{\em Proof.} 
For $k \in \Z_+ $ and $n \in \N$,
set $p_{n,k} := F_n(\{k\})$ and $\tp_{n,k} := p_{n,k}$ for
$k \leq n-2$, with $\tp_{n,n-1} := \sum_{k \geq n-1} p_{n,k}$ and
$\tp_{n,k} :=0$ for $k \geq n$.
Then 
$$
\tPr_{n}[A_r] 
=  \sum_k \tp_{n,k} \prod_{i=1}^k \left( 
\frac{n-i-r}{n-i} \right), 
$$
with the product interpreted as unity for $k=0$.
By Fatou's lemma, 
$\tPr_n[A_r] \to 1$ as  
 $n \to \infty$.
By 
(\ref{eqmonot}) we have  $\tPr_{n}[A_r ] \leq \Pr_{n}[A_r] $
so
$\Pr_{n}[A_r] \to 1$ as well, which
gives us (\ref{limAvoid}).  

By the union bound we have
\bean
 \tPr_n[ H_{r,s}] 
\leq  \sum_k \tp_{n,k} ( k s/(n-1)) \leq (n-1)^{-1} s \mu_{F_n} 
\eean
so that
$\limsup (n \tPr_n[H_{r,s} ]) \leq s \mu_\infty $, and
therefore also 
$\limsup (n \Pr_n[H_{r,s} ]) \leq s \mu_\infty $.
Hence
\bea
\limsup (n \Pr_n[H_{r,s} \cap A_r]) \leq s \mu_\infty ;
~~~~~~~~
\limsup (n \tPr_n[H_{r,s} \cap A_r]) \leq s \mu_\infty .
\label{M17}
\eea
%

By conditioning on the outdegree of vertex 1 
and then using the estimate $e^x \geq 1+x$ for $x \in \R$, we have  
\bea
 \Pr_n[H_{r,s} \cap A_r] 
=  \sum_k p_{n,k} \left( 1 - \frac{ r}{n} \right)^k \left[
1 - \left( 1 - \frac{s}{n-r}\right)^k \right]
\nonumber \\
\geq  \sum_k p_{n,k}
 \left( 1 - \frac{ r}{n} \right)^k \left[
1 - \exp( - ks/(n-r))  \right].
\label{sqbr0}
\eea
Also
\bea
\tPr_n [H_{r,s} \cap A_r] = \sum_k \tp_{n,k} \left[
 \prod_{i=1}^k \left( 1 -  \frac{r}{n-i} \right) -
 \prod_{i=1}^k \left( 1 - \frac{r+s}{n-i} \right) \right].
\label{sqbr}
\eea
Suppose $F_n = F$ for all $n$.
Then both in (\ref{sqbr0}) and (\ref{sqbr}),
 the expression inside  the square brackets is asymptotic
to $n^{-1} ks$,
and hence by Fatou's lemma,
$\liminf (n \Pr_n[H_{r,s} \cap A_r]) \geq s \mu_F$ and 
$\liminf (n \tPr_n[H_{r,s} \cap A_r]) \geq s \mu_F$. 
Combined with (\ref{M17}), this gives us
 (\ref{limHit}) in the case with $F_n = F$ for all $n$.

Now suppose $F_n$  varies with $n$ but
 $F_n \in \cM_{b_n}$ for all $n$, with
$b_n = o(n)$. 
By (\ref{sqbr0}),
\bean
n \Pr_n[H_{r,s} \cap A_r] 
\geq (1-r/n)^{b_n} \sum_k np_{n,k } 
[ 1 - e^{-ks/(n-r)} ]
\eean
and by Taylor's theorem, for $k \leq b_n$ we have
for some $\theta = \theta (n,k) \in (0,1)$ that
$$
 1 - e^{-ks/(n-r)}  = \frac{ks}{n-r} - (1/2) \left( \frac{ks }{n-r} \right)^2 
e^{- \theta ks/(n-r)}   
$$
so that
$$
(n/(ks))
 (1 - e^{-ks/(n-r)} )  \geq 
\frac{n}{n-r}
 -   \frac{nks }{(n-r)^2}  
\geq
1 -   \frac{nb_ns }{(n-r)^2}  .
$$
Hence
$$
n \Pr_n[H_{r,s} \cap A_r] 
\geq
(1-r/n)^{b_n} \left( 1 - \frac{nb_ns}{(n-r)^2} \right) \sum_k  ks p_{n,k},
$$
so that
\bea
\liminf (n \Pr_n[H_{r,s} \cap A_r]) \geq s \mu_\infty .
\label{M18}
\eea 

Next we estimate the right hand side of
(\ref{sqbr}).
By Taylor's theorem we have
\bean
\prod_{i=1}^k \left( 1 - \frac{r}{n-i} \right) \geq \left(
1 - \frac{r }{n-k} \right)^k
\geq 1 - \frac{kr}{n-k} 
\eean
and
\bean
\prod_{i=1}^k \left( 1 - \frac{r+s}{n-i} \right) \leq \left(
1 - \frac{r +s }{n} \right)^k
\leq 1 - \frac{k(r+s)}{n} + \frac{ k^2 (r+s)^2}{ 2 n^2} . 
\eean
Combining these estimates gives us
\bea
\frac{n}{k} 
 \left[
 \prod_{i=1}^k \left( 1 -  \frac{r}{n-i} \right) -
 \prod_{i=1}^k \left( 1 - \frac{r+s}{n-i} \right) \right]
\geq
s - \frac{rk}{n-k}  - \frac{(r+s)^2 k}{n}.
\eea
Hence by (\ref{sqbr}), for $F_n \in \cM_{b_n}$ we have
$$
n \tPr_n [H_{r,s} \cap A_r] \geq \left( \sum_k ks \tp_{n,k} \right) (1 + o(1)), 
$$
and hence
$\liminf (n \tPr_n[H_{r,s} \cap A_r]) \geq s \mu_\infty $. 
Combined with (\ref{M18}) and (\ref{M17}) this gives us 
(\ref{limHit}) in the case with $F_n \in \cM_{b_n}$, completing
the proof.
$\qed$ \\

Let $(Z_m)_{m \geq 0} = (Z_m(F))_{m \geq 0}$ and 
 $(Z'_m)_{m \geq 0} = (Z'_m(\mu_\infty))_{m \geq 0}$ be
branching processes as described in Section \ref{secstate}.
We always assume these branching processes
are independent of each other.
For later use,
we set $T= \sum_{m=0}^\infty Z_m$ and
$T' = \sum_{m=1}^\infty 
Z'_m$.
For $m \in \Z_+$  set
 $\bZ_m := (Z_0,\ldots,Z_m)$
 and $\bZ_m := (Z'_0,\ldots,Z'_m)$.
 Let  $(\tilde{\bZ}_m, \tilde{\bZ}'_m)$ denote an independent copy of
   $({\bZ}_m, {\bZ}'_m)$.

\begin{lemma}
\label{lembranch}
Let $m \in \N$. Then as $n \to \infty$
we have
\bea
({\bf S}_{1,m}(G_{n,F_n}),  
{\bf S}_{2,m}(G_{n,F_n}),  
\bS'_{1,m} (G_{n,F_n}) , 
\bS'_{2,m} (G_{n,F_n}) )  
\tod
(\bZ_m, 
 \tilde{\bZ}_m,
\bZ'_m,
 \tilde{\bZ}'_m) .
\label{M13}
\eea
Also, (\ref{M13}) holds with
$G_{n,F_n}$ replaced by $\tG_{n,F_n}$.
\end{lemma}
{\em  Proof.} We give the argument for $G_{n,F_n}$;
the argument for $\tG_{n,F_n}$ is just the same.

It is rather obvious that $(\bS_{1,m}, \bS_{2,m})$ 
converges in distribution to $(\bZ_m,\tilde{\bZ}_m)$.
Formally, this can be proved by induction on $m$, using
(\ref{limAvoid}).

Suppose we are given (for fixed $m$) the values of 
 $(\bS_{1,m}, \bS_{2,m})$ and consider
for $r \in \N$
 the conditional
distribution of 
 $(\bS'_{1,r}, \bS'_{2,r})$.
We need to
show that this converges to the  distribution 
of $(\bZ'_r,\tilde{\bZ}'_r)$. This is done
by induction in $r$ and we consider the inductive
step, so suppose we also fix for some $r$ the
values of 
 $(\bS'_{1,r}, \bS'_{2,r})$.
Then the value of $S'_{1,r+1}$ is the number of
vertices $j \in [n ] \setminus
 \cup_{s \leq r} S'_{1,s}$ such that $j \to  S'_{1,r}$
(where we use  notation $S'_{1,s}$ to mean either a set of vertices and
its cardinality).

Given $(\bS_{1,m}, \bS_{2,m}, \bS'_{1,r},\bS_{2,r})$,
the number of $j \in \cup_{i \leq m} (S_{1,i} \cup S_{2,i})$
is fixed and the (conditional) probability that any of these has 
 $j \to S'_{1,r}$ tends to zero.
We need to consider the other $j$, i.e. with
$j \notin \cup_{r'\leq r} S'_{1,r'}$ and 
$j \notin \cup_{i \leq m} ( S_{1,i} \cup S_{2,i})$.

For these values of $j$ the conditioning means we know
there are no arcs from $j$ to the set $\cup_{r' < r} S'_{1,r'}$,
a  fixed number  of vertices.  
Therefore
by Lemma \ref{lemPoconv}, 
 the conditional
probability that there is an arc from $j$ to one of
the vertices in $S'_{1,r}$ is asymptotic to  $n^{-1} S'_{1,r} \mu_\infty$.
Hence by standard binomial-Poisson convergence the
(conditional) distribution of the number of such vertices $j$ 
such that $j \to S'_{1,r}$ is asymptotically Poisson
with parameter $ S'_{1,r} \mu_\infty $, which is
the same as the conditional
 distribution of the next value of the branching process
$Z'_{r+1}$.

We can then apply a similar  argument
for $S'_{2,r}$ to complete the induction.
 \qed \\

We now prove a part of Theorem \ref{CDSthm}.

\begin{lemma}
\lbl{CDSpart}
Under the assumptions of
Theorem \ref{CDSthm}, the first conclusion
(\ref{CDSeq}) holds
in the case where $\mu_\infty < \infty$.
\end{lemma}
{\em Proof.}
First consider the graphs $G_{n,F_n}$ with $F_n= F$
for all $n$ and some fixed distribution $F$ on $\Z_+$.
In this case, we can obtain (\ref{CDSeq})
from a result from \cite{CDS}.
The model in \cite{CDS}
 is not described there in terms of
a random digraph, but it is not hard to see that
it can be interpreted that way.
In particular,  the  random variable
$N_n(\tau_n)$ in \cite[Theorem 2.2]{CDS}
 can be interpreted as being the same as our $T_{1}(G_{n,F})$.
Therefore by
 \cite[Theorem 2.2]{CDS},
there exists a coupling of the branching process
$(Z_m)_{m \geq 0}$ and the sequence of random digraphs
$(G_{n,F})_{n \geq 0}$
such that 
\bea
n^{-1} T_1(G_{n,F})
\toP \sigma'(\mu_F) 
{\bf 1} \{ T= \infty \} 
\label{CDStoP}
\eea 
as $n \to \infty$.
Note that our $\sigma'(\mu_F)$ is the $p$ of \cite{CDS}.
The distributional convergence (\ref{CDSeq}) is immediate
from (\ref{CDStoP}).

Next, we consider $G_{n,F_n}$ in the case with 
$F_n$ varying with $n$, assuming also
that  $\mu_\infty < \infty$.
The proof for this case involves adapting the proof
  of
\cite[Theorem 2.2]{CDS}.

The argument in \cite{CDS} (for fixed $F$)
involves considering an exploration process
of the random graph starting from vertex 1,
where at each step one of the currently
unassigned arcs out of one of the vertices
currently being considered is assigned its
destination (uniformly at random over $[n]$),
and if this destination is a previously unconsidered
vertex, then this vertex is added to those currently
being considered at the next stage. If there are
no unassigned arcs out of the current set of
vertices under consideration, the exploration process terminates.

Let $(K_{n,i})_{i \in \N}$
(respectively  $(K_{i})_{i \in \N}$)
be a sequence of independent identically distributed
random variables with the distribution $F_n$
(respectively, $F$).
For $t \geq 0$ 
set $R_n(t) := \sum_{i=1}^{\lfloor t \rfloor } K_{n,i}$ 
and  $R(t) := \sum_{i=1}^{\lfloor t \rfloor } K_{i}$, 
as in (25) of \cite{CDS}.

Let $\surv_n$
denote the event that $1 + R_n(t) -t >0$ for all $t \in \N$
(or equivalently, that $R_n(u+1) -u > 0$ for all
$u \in \Z_+$),
and let $\surv$ 
denote the event that $1 + R(t) -t >0$ for all $t \in \N$.
Note that $\Pr[ \surv_n ] = \sigma(F_n)$ 
and $\Pr[ \surv ] = \sigma(F)$, because the
exploration process of a branching process
with offspring distribution $F$
can be interpreted as a random walk with 
successive steps having the distribution of $K_1 -1$. 

For $t \in \Z_+$, let $N_n(t)$ denote the number
of coupons collected after $n$ attempts in a coupon
collector process with $n$ coupons (starting with $N_n(0)=1$;
for a formal description see \cite{CDS}),
running independently of the random walk $R_n(\cdot)$.
For $ t \geq 0$, set $S_n(t) :=  R_n( N_n(\lfloor t \rfloor)  ) - t$.
As described in in \cite{CDS}, there is a coupling in which 
$S_n(t)$ (for $t \in \N$) can
 be viewed as  the total number of unassigned out-arcs
from the current set of vertices after $ t $
stages of the exploration process,  up to 
time $\tau_n$, where $\tau_n$ denotes
the first  $t $ such that $S_n(t) \leq 0$.

We claim that
there exists $\eps >0$ such that
\bea
\lim_{ n \to \infty} (
\Pr [ \tau_n \geq n\eps, \surv_n ] ) = \sigma(F).
\label{M30}
\eea
 This is proved by following the proof of
Lemma 4.2 of \cite{CDS}
(the notation $\sigma^{\rm{GW}}$  in \cite{CDS}
denotes an extinction  probability, whereas our $\sigma(F)$ is
a survival probability!)
Most of the proof of Lemma 4.2 of \cite{CDS} carries over
easily to the present setting. We just elaborate
on the assertion in that proof that
$z^{-1}G(z)^{1-2 \eps} <1$ for some $z <1$. Here the
$G$ of \cite{CDS} is a probability generating function which we
denote by  $\phi_n$, with $\phi_n(z) := \sum_k z^k F_n(\{k\})$.  
Also set $\phi(z) := \sum_k z^k F(\{k\})$, and  
note that $\lim_{n \to \infty} \phi_n(z) = \phi(z) $ for
$z \in (0,1)$.

By Fatou's lemma,
$
\liminf_{a \downarrow 0} a^{-1}(1- \phi(1-a)) \geq \mu_F.
$ 
Assuming $\mu_F >1$, taking $\eps >0$ and
$\delta >0$ with $(1-2 \eps) (\mu_F - \delta) >1$,
and then $a \in (0,1)$ (close to 0) with
$(1- a (\mu_F-\delta))^{1-2 \eps} < (1-a)$
and also
$a^{-1} ( 1- \phi(1-a) ) >  \mu_F - \delta$,
we have that
$\phi(1-a)^{1-2 \eps} < (1-a)$, and hence
for large $n$ that $ \phi_n(1-a)^{1-2 \eps} < 1-a$.

For $q >0$ we have $n^{-1} \E[R_n(nq) ] \to q \mu_\infty$ as
$n \to \infty$, and
since we are assuming $\mu_\infty < \infty$ and $b_n = o(n)$, we have
\begin{align*}
\Var \left[ n^{-1} R_n(nq) \right] =  n^{-2} \lfloor nq \rfloor 
\Var[ K_{n,1}] 
& \leq (q/n) \E K_{n,1}^2
\\
&  \leq (q/n) \E [ b_n K_{n,1} ] 
\to 0.
\end{align*}
By following
the proof (26) of \cite{CDS}, 
for each positive $s$ we have that 
\bea
n^{-1} S_n (ns) \to (1- e^{-s}) \mu_\infty -s
\label{M20}
\eea   
in probability.
This weaker version of (26) of \cite{CDS} suffices to
give us (27) of \cite{CDS}.

At the end of the 3-line display just 
after (27) of \cite{CDS}, there are three
terms which we wish to show tend to zero.
The first term tends to zero by (\ref{M30}) and
the fact that $\sigma(F_n) \to \sigma(F)$ as $n \to \infty$.
The third term can be shown to tend to zero using
the same fact.
To show the second term tends to zero, we use the
next three-line display of \cite{CDS}; we need to
check that for $\delta >0$ we have
\bea
\limsup_{n \to \infty} \Pr [ \inf S_n(ns), s \in [\eps, \theta - \delta ] ]
=0,
\label{M21}
\eea
where $\theta$ is the solution in $(0,\infty)$ to
$(1-e^{-\theta})/\theta = 1/\mu_\infty$.
To see this, set
$h: = \inf_{s \in [\eps, \theta - \delta]} ((1- e^{-s} ) \mu_\infty -s) >0$
and take 
$s_1,\ldots,s_\ell \in [\eps,\theta-\delta]$ with $s_1 = \eps$,
$s_\ell = \theta -\delta$ and $0 < s_{i+1} -s_i < h/4$
for $1 \leq i \leq \ell -1$. By (\ref{M20}),
with probability tending to 1 we have for 
for each $i$ that $n^{-1} S_n(n s_i) > h/2$, and then
using that $S_n(t) = R_n(N_n(\lfloor t \rfloor) -t)$ 
we have $n^{-1} S_n(ns) \geq h/4$ for all $s \in [s_i,s_{i+1}]$
which gives us $S_n(ns) >0$ for $s \in [\eps,\theta -\delta]$, and
hence (\ref{M21}).  We can then follow the rest
of the argument in \cite{CDS} to get (\ref{CDSeq}) in the case
where $\mu_\infty < \infty$.  
$\qed$

\begin{lemma}
\label{CDSlem}
It is the case that 
\bea
\lim_{n \to \infty}
\Pr_n[1 \tovia 2]
=
  \sigma'(\mu_\infty) \sigma(F)
\label{M26}
\eea
and
\bea
\lim_{n \to \infty}
\tPr_n[1 \tovia 2] =
  \sigma'(\mu_\infty) \sigma(F).
\label{M27}
\eea
\end{lemma}
{\em  Proof.}
Write just $\sigma'$ for $\sigma'(\mu_\infty)$ and
$\sigma$ for $\sigma(F)$,
 and $T_1$ for $T_1(G)$.
We first prove (\ref{M26}).
By symmetry we have $\Pr_n[1 \tovia 2 | T_1] = (T_1-1)/n$, and therefore
by Lemma \ref{CDSpart}, we have
\bea
\lim_{n \to \infty} \Pr_n[1 \tovia 2] = \lim_{n \to \infty}
\E_n [ ( T_1 -1)/n]  = \sigma'(\mu_\infty) \sigma(F), ~~~~
\mbox{ if } \mu_\infty < \infty. 
\label{M11}
\eea

Now suppose $\mu_\infty = \infty$. Given $\eps >0$ we can choose
$K \in \N$ such that $\Pr[ T >K] < \sigma + \eps$.
By branching process approximation (Lemma \ref{lembranch}),
$
\lim_{n \to  \infty} \Pr_n [  T_1 > K 
] = \Pr [  T > K   ] < \sigma + \eps,
$
and also by symmetry 
$\Pr_n [ \{ 1 \tovia 2 \} \cap \{T_1 \leq K \} ] \leq
 (K/(n-1)) \Pr_n[ T_1 \leq K ] $,
which tends to zero,
so
\bea
\limsup_{n \to \infty}  \Pr_n[1 \tovia 2 ] \leq \sigma  + \eps.
\label{M25}
\eea 
 Given $h \in \N$, let
$F^h$ (respectively $F_n^h$) denote the distribution  of a random variable
$\min(\xi , h)$ (respectively $\min(\xi^n,h)$), where 
$\xi$ (respectively $\xi^n$) is a random variable with
distribution $F$ (respectively $F_n$). Pick $h \in \N$ with
$\sigma'(h) > 1-\eps$ and $ \sigma(F^h) \geq \sigma(F)(1 - \eps)$.
Here we are using the continuity of the branching process
survival probability in the offspring distribution.

Given $n$, choose $a_n \in \N$ with
$\mu_{F_n^{a_n}}  \in [h,h+1] $ (this is possible
for all large enough $n$ because $\mu_\infty = \infty$).
Note that $a_n \geq h$.
Let $\Pr_n^*$ denote probability for a random digraph of the
form of
$G_{n, F_n^{a_n}}$.

Suppose first that $a_n \to \infty$ as $n \to \infty$. Then $F_n^{a_n}$
converges weakly to $F$, so
 by monotonicity and the case already proved,
we have
$$
\liminf \Pr_n[ 1 \tovia 2 ] \geq \liminf \Pr^*_n [ 1 \tovia 2]
\geq \sigma'(h) \sigma(F) \geq (1- \eps) \sigma(F).
$$

Suppose instead that $a_n$ is bounded. For
any subsequence of $n$ we can take a further subsequence,
such that 
along this subsequence 
$a_n$ tends to a finite limit $a$ 
so that
$F_n^{a_n}$
converges weakly to $F^a$, and also $\mu_{F_n^{a_n} } $ tends
to a limit $y$ (between $h$ and $h+1$).
Also $a \geq h$ so $\sigma(F^a) \geq \sigma(F^h) \geq \sigma(F) (1-\eps)$.
Then by monotonicity and the case already proved,
as $n \to \infty$
along this further subsequence we have
$$
\liminf \Pr_n[ 1 \tovia 2 ] \geq \liminf \Pr^*_n [ 1 \tovia 2]
=  \sigma' (y) \sigma (F^{a}) \geq (1-\eps)^2 \sigma(F),
$$
and since $\eps$ is arbitrary, combined with (\ref{M25}) this gives us
(\ref{M26}) for the case $\mu_\infty = \infty$.
Combined with (\ref{M11}) this gives us (\ref{M26}) in full generality.

Now consider $\tG_{n,F_n}$. 
By (\ref{eqmonot}) 
and
 (\ref{M26}) we have
\bea
\liminf_{n \to \infty} \tPr_n [1 \tovia 2] \geq \liminf_{n \to \infty} \Pr_n [ 1 \tovia 2]
= \sigma'(\mu_\infty) \sigma(F).
\label{M12}
\eea
On the other hand, given $\eps >0$
 we can choose $K \in \N$ such that
$\Pr [ \{ T > K \}  \cap \{ T' >K \} ]   <  \sigma'(\mu_\infty) \sigma(F) + \eps $. 
By branching process approximation (Lemma \ref{lembranch}),
$$
\lim_{n \to  \infty} \tPr_n [ \{ T_1 \leq K \} \cup \{ T'_2 \leq K \}
] = \Pr [ \{ T \leq K  \} \cup \{T' \leq K \}  ] > 1- (\sigma' \sigma + \eps),
$$
and also by symmetry 
$\tPr_n [ \{ 1 \tovia 2 \} \cap \{T_1 \leq K \} ] \leq
 (K/(n-1)) \tPr_n[ T_1 \leq K ] $ which tends to zero,
and similarly
$\tPr_n [ \{ 1 \tovia 2 \} \cap \{T'_2 \leq K \} ] \to 0$,
so by the union bound
$$
\lim_{n \to \infty} \tPr_n[  \{1 \tovia 2 \}  \cap (
 \{ T_1 \leq K \} \cup \{ T'_2 \leq K \} )] = 0.
$$ 
Therefore 
$\limsup_{n \to \infty}  \tPr_n[1 \tovia 2 ] \leq \sigma' \sigma + \eps$. 
Combined with (\ref{M12}) this shows that
$\lim_{n \to \infty}  \tPr_n[1 \tovia 2 ] = \sigma' \sigma $. 
\qed \\

{\em Proof of Theorem \ref{CDSthm}.} 
Set $X_n = n^{-1} T_1(G_{n,F_n})$ and
$\tilde{X}_n := n^{-1} T_1(\tG_{n,F_n})$. 
Given $\eps >0$,
we may choose finite $K$ such that $\Pr[ T >K ] \leq \sigma(F) + \eps/2$.
Then by branching process approximation (Lemma \ref{lembranch}),
 for large enough $n$ we have
$
\Pr [ \tilde{X}_n \leq \eps ]  \geq \Pr[ 
T_1(\tG_{n,F_n}) \leq K ] \geq 1 - \sigma(F)  - \eps.
$
Also 
 $\tilde{X}_n$ stochastically dominates
$X_n$ by (\ref{eqmonot}).
Hence 
\bea
\liminf_{n \to \infty} \Pr[ X_n \leq t ]
\geq
\liminf_{n \to \infty} \Pr[ \tilde{X}_n \leq t ]
  \geq 1 - \sigma(F), ~~~~ t >0.
\label{M22}
\eea

In view of Lemma \ref{CDSpart}, to prove (\ref{CDSeq})
we only need to consider the case
with $F_n$ varying with $n$ and $\mu_\infty = \infty$,
so we assume $\mu_\infty = \infty$ for  a while.
 Then $\sigma'(\mu_\infty)=1$.
 By Lemma  \ref{CDSlem}
we have 
\bea
\lim_{n \to \infty} \E X_n = \sigma (F).
\label{M23}
\eea
In the case $\sigma(F) = 0$ this gives us (\ref{CDSeq}) at once,
so now assume $\sigma(F) >0$ too.
Let $\eps \in (0,\sigma(F))$.
Since $X_n \leq 1$ we have $\E X_n \leq \eps \Pr[X_n \leq \eps]
+ (1- \Pr[X_n \leq \eps])$, so that
$
(1- \eps) \Pr[X_n \leq \eps] \leq 1 - \E X_n
$
and using (\ref{M23}) we have
\bean
\limsup_{n \to \infty} \Pr[ X_n \leq \eps] \leq (1-\eps)^{-1} (1- \sigma(F)) 
\eean
so that
\bea
\liminf_{n \to \infty} \Pr[ X_n > \eps] \geq (1-\eps)^{-1} (\sigma(F) - \eps)
\geq \sigma(F) - \eps . 
\label{M24}
\eea
Let $\delta \in (0,1/2)$ and set $\eps = \sigma(F) \delta^2/2$. 
Suppose $\Pr [ X_n \leq 1- \delta | X_n > \eps ] > \delta$ for infinitely
many $n$. Then for such $n$ we have
$$
\E[ X_n| X_n > \eps] \leq (1- \delta) \delta + (1-\delta) = 1 - \delta^2
$$
and hence by (\ref{M22}),  along this subsequence
\bean
\limsup_{n \to \infty} \E X_n = \limsup \left( 
\E[X_n {\bf 1}\{X_n \leq \eps\} ] + \Pr[X_n > \eps] \E[X_n| X_n > \eps ]
\right)
\\
\leq \eps  + \sigma(F) (1 - \delta^2) = \sigma(F) ( 1 - \delta^2/2)   
\eean
which contradicts (\ref{M23}).  Hence $\Pr[X_n \leq 1-\delta| X_n > \eps] 
\leq \delta$ for all but finitely many $n$, and using (\ref{M24})
we have
\bean
\liminf_{n \to \infty} \Pr[ X_n > 1- \delta] \geq
(1 - \delta)
 (1-\delta^2/2)
\sigma(F).
\eean 
Hence for $t \in (0,1)$,
we have $
\liminf_{n \to \infty} \Pr [ X_n > t] \geq  \sigma (F)$,
and with (\ref{M22}) this shows
that $\Pr[X_n \leq t] \to 1-\sigma (F)$.
This gives us (\ref{CDSeq}).

We still need to prove (\ref{CDS2eq}), not only under the assumption
$\mu_\infty = \infty$ so now relax this assumption.
By (\ref{CDSeq}),
for $ t < \sigma'(\mu_\infty) $
we have
$$
\limsup \Pr [ \tilde{X}_n \leq t] \leq \limsup \Pr[X_n \leq t] 
\leq 1 - \sigma(F).
$$
Hence by (\ref{M22}), 
\bea
\lim_{n \to \infty} \Pr[ \tilde{X}_n \leq t] = 1 - \sigma(F), ~~~ \mbox{ for } 
0 < t <
 \sigma'(\mu_\infty) .
\label{M29}
\eea

Next, let $T' := \sum_{m \geq 0} Z'_m(\mu_\infty)$ as before.
 Given $\eps >0$, choose $K $ with 
$\Pr[T' > K ] \leq \sigma'(\mu_\infty) + \eps/2$. 
 let $N_{{\rm small}} = \sum_{i=1}^n {\bf 1}
 \{ T'_i(\tG_{n,F_n}) \leq K \} $.
Using Lemma \ref{lembranch} we have that
 $$
\tPr_{n,F_n}[ T'_1
\leq K]
\to \Pr [T' \leq K]; ~~~~
 \tPr_{n,F_n} [ T'_1
\leq K
 , T'_2
 \leq K]
\to (\Pr [T' \leq K])^2,
$$
 and hence
$\E[ n^{-1} N_{{\rm small}} ] \to \Pr[T' \leq K] $
and
$\Var[ n^{-1} N_{{\rm small}} ] \to 0 $.
Hence 
$$
 \Pr[ N_{{\rm small}} \leq n(1 - \sigma'(\mu_\infty) -  \eps ) ] \leq
\Pr[ N_{{\rm small}} \leq n(\Pr[T' \leq K] - \eps/2 ) ] \to 0.
$$
Given $n$,
let $I_{{\rm small}}$ be the set of indices $j \in [n]$
such that
 $ T'_j(\tG_{n,F_n}) \leq K  $.
Then 
$$
\sum_{i=1}^n {\bf 1} 
\left\{ 
\sum_{j \in I_{{\rm small}} } {\bf 1}
 \{ i \tovia j \mbox{ in } \tG_{n,F_n} \} > \eps n
 \right\} \leq K/\eps,
$$
so that by symmetry
$$
\tPr_{n,F_n} 
\left[ 
\sum_{j \in I_{{\rm small}} } {\bf 1}
 \{ 1 \tovia j\} > \eps n
 \right] \leq K/(n \eps),
$$
and hence setting $\sigma' = \sigma'(\mu_\infty)$, by the union bound
we have
\bean
\Pr [ \tilde{X}_n \geq \sigma' + 2 \eps ]
 \leq \Pr[ N_{{\rm small}} \leq n ( 1 - \sigma' -\eps) ]
+  \tPr_{n, F_n}
\left[ 
\sum_{j \in I_{{\rm small}} } {\bf 1}
 \{ 1 \tovia j\} > \eps n
 \right] \to 0. 
\eean
Combined with (\ref{M29}) this gives us (\ref{CDS2eq}).
\qed

\section{Proof of Theorems \ref{thgiant}, \ref{thFvary} and \ref{thunif}}

In this section, 
we make the same assumptions about $F,F_n, \mu_\infty, $ and
(if applicable) $b_n$, and use the same notation
$\Pr_n$, $\tPr_n$, $\E_n$, $\tE_n$, as we did in
the previous section.
Also $T$ and $T'$ are as in the previous section,
and we set $\sigma:= \sigma(F)$ and $\sigma':=\sigma'(\mu_\infty)$.

\begin{lemma}
\label{L1}
It is the case that 
\bea
\lim_{n \to \infty} \P_n[1 \leftrightarrow 2 ] = (\sigma' \sigma)^2;
\label{M1}
\\
\lim_{n \to \infty} \tPr_n[1 \leftrightarrow 2 ] = (\sigma' \sigma)^2.
\label{M1a}
\eea
\end{lemma}
{\em Proof.}
We prove (\ref{M1}); the proof of (\ref{M1a}) is just the same
but with $\Pr_n$ replaced by $\tPr_n$ throughout.

If $\mu_F  \leq 1$ then $\sigma' \sigma=0$ and (\ref{M1}) follows from
Lemma \ref{CDSlem}, so now assume $\mu_F  >1$.
Then
 $\mu_\infty \geq \mu_F > 1$
 by Fatou's lemma,
so $\sigma' \sigma >0$. 
Choose $K$ such that $\Pr[ \{ T > K \}  \cap \{T' > K \} ] 
\approx \sigma' \sigma$.
The $\approx$ means the approximation can be made arbitrarily
close to zero by the choice of $K$ (but we now view $K$ as fixed).

Let $T_i,T'_i$ be as in Section 
\ref{secstate}. Then
by symmetry,
$
 \Pr_n[ \{ 1 \tovia 2 \} \cap \{ T_{1}  \leq K \} ]  \leq
(K/(n-1)) \Pr_n[ T_1 \leq K],
$
which tends to zero as $n  \to \infty$.  Similarly
$
\lim_{n  \to \infty} \Pr_n[ \{ 1 \tovia 2 \} \cap \{  T'_{2}  \leq K \} ]  = 0 
$
so that
\bea
\lim_{n  \to \infty} \Pr_n[ \{ 1 \tovia 2 \} \setminus (
 \{ T_{1}  > K \}  \cap \{ T'_{2}  >  K \} ) ] = 0. 
\label{M10}
\eea
Also, by branching process approximation (Lemma \ref{lembranch}), 
$$
\lim_{n \to \infty} \Pr_n [ \{ T_{1} >  K \}
\cap \{ T'_{2} > K \} ]
= \Pr[ \{ T > K \} \cap \{T' > K\} ]
\approx \sigma' \sigma.
$$ 
Therefore by (\ref{M10}) and Lemma \ref{CDSlem},
with $\triangle $ denoting the symmetric difference of two sets,
$$
\limsup_{n  \to \infty} \Pr_n[ \{ 1 \tovia 2 \} \triangle (
 \{ T_{1}  > K \}  \cap \{ T'_{2}  >  K \} ) ]  \approx 0. 
$$
Similarly,
$$
\limsup_{n  \to \infty} \Pr_n[ \{ 2 \tovia 1 \} \triangle (
 \{ T_{2}  > K \}  \cap \{ T'_{1}  >  K \}
 ) ]  \approx 0, 
$$
and therefore
$$
\limsup_{n  \to \infty} \Pr_n[ \{ 1 \leftrightarrow 2 \} \triangle (
 \{ T_{1}  > K \}  \cap \{ T'_{1}  >  K \}
\cap
 \{ T_{2}  > K \}  \cap \{ T'_{2}  >  K \} ) ]  \approx 0. 
$$
By 
 branching process approximation (Lemma \ref{lembranch}),
and the inclusion-exclusion formula,
\bean
\lim_{n \to \infty}
 \Pr [ \{ T_{1}  > K \}  \cap \{ T'_{1}  >  K \}
\cap
 \{ T_{2}  > K \}  \cap \{ T'_{2}  >  K \}  ]
\\
=
 (\Pr [ T  > K   ])^2 ( \Pr  [ T'  >  K ])^2
\approx (\sigma' \sigma)^2,
\eean
and the result (\ref{M1}) follows. $\qed$ \\

{\em Proof of Theorems \ref{thgiant} and \ref{thFvary}.}
We simultaneously prove (\ref{M6}) and (\ref{M15});
the proof of (\ref{L2conv}) and (\ref{M15a})
 is just the same (with $\Pr_n$ replaced by $\tPr_n$ throughout).
Choose a large constant $K$ 
such that
\bean
\Pr [ T > K ] \approx  \sigma; ~~~~~~~~ 
\Pr [T' > K ] \approx  \sigma', 
\eean
with $\approx$ interpreted as in the preceding proof.
Given $n$, define the events
 \bean
E_i :=
\{ T_i \leq K \} \cup \{ T'_i  \leq K \}, ~~~~ i \in [n].
\eean
By branching process approximation (Lemma \ref{lembranch}), as $n \to \infty$ we have
\bean
\Pr_n[ E_1] \to
 \Pr[ \{ T \leq K\} 
\cup
   \{ T' \leq K 
\}
  ];
\\
\Pr_n[ E_1 \cap  E_2] \to
 ( \Pr[  \{ T  \leq K \}  \cup \{T' \leq K \}  
 ] )^2 .
\eean
 Therefore, setting
$ N_{\rm small} := \sum_{i=1}^n {\bf 1}_{E_i},  $
we have
$\Var(N_{\rm small}/n) \to 0$ and
\bea
n^{-1} N_{\rm small} \toP \Pr [ 
\{ T \leq K \} \cup \{ T' \leq K \} ] \approx
1 - \sigma' \sigma.
\label{M7}
\eea
Suppose first that $\mu_F \leq 1$. Then $\sigma' \sigma =0$,
and given any $\eps >0$ we may choose $K$ such that  if
 $n > K/\eps$ we have
 by (\ref{M7}) that
$$
\Pr_n [ L_1(G) > \eps n] \leq \Pr_n [ N_{\rm small} <  (1 - \eps) n] 
\to 0.
$$
This gives us
(\ref{M6}) and
(\ref{M15})
 in the case where $\mu_F  \leq 1$.

Now suppose $ \mu_F > 1$. Then $\sigma = \sigma(F)>0$,
and by Fatou's inequality $\mu_\infty \geq \mu_F >1$ so
$\sigma' = \sigma'(\mu_\infty) > 0$.
Let $N_{>K}$  be the number of vertices of $G_n$ lying in
strongly connected components of order greater than $K$.
 For $i \geq 1$,
let us write just $L_i$ for $L_i(G)$.
Let $I$ be the last  $i$ such that $L_i > K$.
Then
 \bea
\Pr_n[ 1 \leftrightarrow 2 | (L_1,L_2,\ldots)] \leq
\sum_i \left(\frac{L_i}{n} \right)^2 
\leq \sum_{i \leq I} \left(\frac{L_i}{n}\right)^2
+ \sum_{I <i \leq n} \left( \frac{K}{n}\right)^2
\nonumber \\
\leq \frac{L_1}{n} \left( \sum_{i \leq I} \frac{L_i}{n} \right)
+ \frac{K^2}{n} 
\nonumber \\
= \frac{L_1}{n} \times \frac{N_{>K}}{n} + \frac{K^2}{n}.
\label{M4}
\eea
Note that $N_{>K}$ is determined by $(L_1,L_2,\ldots,)$.
Let $\eps \in (0,1)$.
By (\ref{M4}) we have
\bea
\Pr_n[1 \leftrightarrow 2 | N_{>K} < (1+\eps^2) \sigma' \sigma n,
L_1 \leq (1-\eps) \sigma' \sigma n] \leq(1-\eps) (1+\eps^2) (\sigma' \sigma)^2
+ o(1)
\label{M8}
\eea
and (using $L_1 \leq \max(N_{>K},K)$) also
\bea
\Pr_n[1 \leftrightarrow 2 | N_{>K} < (1+\eps^2) \sigma' \sigma n,
L_1 > (1-\eps) \sigma' \sigma n] \leq (1+\eps^2)^2 (\sigma' \sigma)^2
+ o(1).
\label{M9}
\eea

Now $N_{>K} \leq n - N_{\rm small}$   so 
by (\ref{M7}), given $\eps >0$ we can choose $K$ 
so that
\bea
\Pr_n [ n^{-1} N_{>K} < (1+ \eps^2) \sigma' \sigma ] \to 1.
\label{M5}
\eea
Then by (\ref{M8}) and (\ref{M9})
we have
\begin{align*}
\Pr_n [1 \leftrightarrow 2] \leq 
 (\sigma' \sigma)^2 &
( (1-\eps) (1+ \eps^2) \Pr_n[L_1 \leq (1-\eps) \sigma' \sigma n]
\\
 & + (1+\eps^2)^2 ( 1 - 
\Pr_n[L_1 \leq (1-\eps) \sigma' \sigma n] ) ) 
+o(1)
\\
\leq (\sigma' \sigma)^2
& ( 1 + \eps ( 4 \eps - \Pr_n [L_1 \leq (1-\eps) \sigma' \sigma n] ) ) +o(1)  
\end{align*}
and by comparison with (\ref{M1}) this shows that
$$
\limsup_{n \to \infty}
 \Pr_n [L_1 \leq (1-\eps) \sigma' \sigma n]  \leq  4 \eps.
$$
Together with (\ref{M5}) and the fact that
$L_1 \leq \max( N_{>K},K)$, this gives us first part of
 (\ref{M6}) and of (\ref{M15}).

Since $L_2 \leq \max
(N_{>K} - L_1,K)$, for $n > K/ (2 \eps^2 \sigma' \sigma)$ 
with $\eps,K$ as in (\ref{M5})  we
have
\bean
\Pr_n [ n^{-1} L_2 > 2 \eps^2 \sigma' \sigma]
\leq \Pr_n [ n^{-1} (N_{>K} - L_1 )  > 2 \eps^2 \sigma' \sigma]
\\
\leq \Pr_n [ n^{-1} N_{>K}  > (1 + \eps^2) \sigma' \sigma]
+ 
 \Pr_n [ n^{-1}   L_1   < (1 - \eps^2) \sigma' \sigma],
\eean
which tends to zero
 by (\ref{M5}) and 
 the first part of (\ref{M6}) or (\ref{M15}).
 This shows that
$L_2/n \toP 0$, which is the second part of (\ref{M6}) and of (\ref{M15}).
  $\qed$ \\

{\em Proof of Theorem \ref{thunif}.}
Suppose (\ref{M16}) fails.
Then we can find a sequence of distributions $F_n \in \cM_{b_n}$
such that 
$\limsup_{n \to \infty} \sup_{F \in \cM_{b_n}}
  d (n^{-1} L_1(G_{n,F_n}) ,    \sigma'(\mu_{F_n}) \sigma(F_n) )  >0$.
By taking  a subsequence, we may assume $F_n$ converges
to a 
limiting distribution $F$ on $\Z_+ \cup \{\infty\}$
and $\mu_{F_n}$ converges to a (possibly infinite)
limit $\mu_\infty$. But then we would have a contradiction
of Theorem \ref{thFvary}.

This gives us  (\ref{M16}). The proof of (\ref{M16a}),
and of the stated results for $L_2(G_n)$ and $L_2(\tG_n)$,
is similar. $\qed$ \\ 

{\em Acknowledgement.} I thank Dieter Mitsche for organising
a pleasant conference in Nice in 2014, where I learnt about
the work in \cite{CDS}.



\begin{thebibliography}{}

\bibitem{BGJ}
Bloznelis, M., G\"otze, F. and Jaworski, J.
 (2012) 
Birth of a strongly connected giant in an inhomogeneous random digraph.
{\em J. Appl. Probab.} {\bf  49},
 601-611. 


\bibitem{Boll}
Bollob\'as, B. (2001)
{\em Random Graphs (Second edition)}.
 Cambridge University Press, Cambridge.

\bibitem{Broder}
     Broder, A.,  Kumar, A. Maghoul, F.,
     Raghavan, P.,  Rajagopalan, S., 
 Stata, R.  Tomkins, A. and Wiener, J.
(2000)
Graph structure in the Web.
{\em Computer Networks} {\bf 33}, 309-320


\bibitem{CDS}
Comets, F.,  Delarue, F. and Schott, R.
(2014)
Information transmission under random emission constraints.
 {\em Combin. Probab. Comput.} {\bf 23}, 973-1009.

\bibitem{CF}
Cooper, C. and Frieze, A.
 (2004)
The size of the largest strongly connected component of a random digraph 
with a given degree sequence. 
{\em Combin. Probab. Comput.} {\bf 13},
 319-337. 

\bibitem{DMS}
Dorogovtsev, S.N.  Mendes, J.F.F. and Samukhin, A.N.
(2001)
Giant strongly connected component of directed networks
{\em Phys. Rev.} E {\bf 64}, 025101(R), 4pp 

\bibitem{Feller}
Feller, W. (1968)
{\em An Introduction to Probability Theory and its Applications. Vol. I.}
3rd ed. John Wiley and Sons,  New York

\bibitem{FF}
Fenner, T. I. and Frieze, A. M. (1982)
On the connectivity of random m-orientable graphs and digraphs.
{\em Combinatorica} {\bf 2}, 
347-359. 


\bibitem{NSW}
 Newman, M.E.J.
 Strogatz, S.H.
 and Watts, D.J. (2001)
Random graphs with arbitrary degree distributions and their applications.
{\em Phys. Rev. E} {\bf 64}, 026118, 17pp.
\end{thebibliography}
\end{document}